\documentclass[12pt]{amsart} 
\usepackage{amssymb,amsmath,amscd} 
\theoremstyle{plain} 
\newtheorem{Lem}{Lemma}[section] 
\newtheorem{Prop}[Lem]{Proposition} 
\newtheorem{Thm}[Lem]{Theorem} 
\newtheorem{Cor}[Lem]{Corollary} 
\theoremstyle{definition} 
\newtheorem{Def}[Lem]{Definition} 
\newtheorem{Rem}[Lem]{Remark} 

\errorcontextlines=0  \numberwithin{equation}{section}

\newcommand{\bbC}{{\mathbb C}} 
\newcommand{\bbR}{{\mathbb R}} 
 
\newcommand{\bbZ}{{\mathbb Z}} 
\newcommand{\bbN}{{\mathbb N}} 
\newcommand{\bbK}{{\mathbb K}} 
\newcommand{\calA}{{\mathcal A}} 
 
\newcommand{\calH}{{\mathcal H}} 
 
\newcommand{\calM}{{\mathcal M}} 
 
\newcommand{\calP}{{\mathcal P}} 
\newcommand{\calS}{{\mathcal S}} 
\newcommand{\calT}{{\mathcal T}} 
\newcommand{\calV}{{\mathcal V}} 
\newcommand{\Id}{\mathrm{Id}} 
\newcommand{\frakl}{{\mathfrak l}} 
\newcommand{\frakr}{{\mathfrak r}}

\newcommand{\bpi}{\begin{picture}} 
\newcommand{\epi}{\end{picture}}

\begin{document} 

\title[Algebraic curvature tensors] 
{Determination of the structure of algebraic curvature tensors by means of Young symmetrizers} 
\date{June 2001} 
\author[B. Fiedler]{Bernd Fiedler}
\address{Bernd Fiedler \\ Mathematisches Institut \\ Universit\"at Leipzig\\ 
Augustusplatz 10/11 \\ D-04109 Leipzig \\ Germany}
\urladdr{http://home.t-online.de/home/Bernd.Fiedler.RoschStr.Leipzig/}  
\email{Bernd.Fiedler.RoschStr.Leipzig@t-online.de}  

\begin{abstract}
For a positive definite fundamental tensor all known examples of Osserman algebraic curvature tensors have a typical structure. They can be produced from a metric tensor and a finite set of skew-symmetric matrices which fulfil Clifford commutation relations. We show by means of Young symmetrizers and a theorem of S. A. Fulling, R. C. King, B. G. Wybourne and C. J. Cummins that every algebraic curvature tensor has a structure which is very similar to that of the above Osserman curvature tensors. We verify our results by means of the Littlewood-Richardson rule and plethysms. For certain symbolic calculations we used the Mathematica packages {\sf MathTensor}, {\sf Ricci} and {\sf PERMS}.
\end{abstract}

\maketitle 

%First page headline in LaTeX for S\'eminaire Lotharingien de Combinatoire
%--first part
%\thispagestyle{myheadings}
%\font\rms=cmr8 
%\font\its=cmti8 
%\font\bfs=cmbx8
%\markright{\its S\'eminaire Lotharingien de
%Combinatoire \bfs 46 \rms (2001), Article~B46xxx\hfill}
%\def\thepage{}
%
%

\section{Introduction}

Let ${\calT}_r V$ be the vector space of the $r$-times covariant tensors $T$ over a finite-dimensional $\bbK$-vector space $V$, $\bbK = \bbR$ or  $\bbK = \bbC$.
\begin{Def}
A tensor $T \in {\calT}_4 V$ is called an {\itshape algebraic curvature tensor} iff $T$ has the {\itshape index commutation symmetry}
\begin{eqnarray}
\forall\,u, x, y, z\in V:\;\;\;
T(u,x,y,z) & = & - T(u,x,z,y) \;=\; T(y,z,u,x)
%\label{riemsym}
\end{eqnarray}
and fulfills the {\itshape first Bianchi identity}
\begin{eqnarray}
\forall\,u, x, y, z\in V:\;\;\;T(u,x,y,z) + T(u,y,z,x) + T(u,z,x,y) & = & 0\,.
\end{eqnarray}
\end{Def}
If we consider the coordinates $T_{i j k l}$ of a tensor $T\in {\calT}_r V$ then $T$ is an algebraic curvature tensor iff its coordinates satisfy
\begin{eqnarray}
T_{i j k l} & = & - T_{j i k l} \;=\; - T_{i j l k} \;=\; T_{k l i j}
\label{riemsym}
\end{eqnarray}
and
\begin{eqnarray}
T_{i j k l} + T_{i k l j} + T_{i l j k} & = & 0 \,.
\end{eqnarray}
with respect to every basis of $V$.

We assume that $V$ possesses a fundamental tensor $g \in {\calT}_2 V$ (of arbitrary signature) which can be used for raising and lowering of tensor indices.
\begin{Def}
Let $T \in {\calT}_4 V$ be an algebraic curvature tensor and $x \in V$ be a vector with $|g(x,x)| = 1$.
The {\itshape Jacobi operator} $J_T(x)$ of $T$ and $x$ is the linear operator
$J_T(x) : V \rightarrow V\;,\; J_T(x): y \mapsto J_T(x)y$ that is defined by
\begin{eqnarray}
\forall\, w \in V : \;\;\;
g(J_T(x) y , w) & = & T(y, x, x, w)\,. \label{jacobi}
\end{eqnarray} 
\end{Def}
Furthermore we use the following operators $\gamma$ and $\alpha$:
\begin{Def}
\begin{enumerate}
\item{Let $S \in {\calT}_2 V$ be a symmetric tensor of order 2, i.e. the coordinates of $S$ satisfy $S_{i j} = S_{j i}$. We define a tensor $\gamma (S) \in {\calT}_4 V$ by
\begin{eqnarray}
\gamma (S)_{i j k l} & := &
{\textstyle \frac{1}{3}} \left( S_{i l} S_{j k} - S_{i k} S_{j l} \right) \,.
\label{gammas}
\end{eqnarray}
}
\item{Let $A \in {\calT}_2 V$ be a skew-symmetric tensor of order 2, i.e. the coordinates of $A$ satisfy $A_{j i} = - A_{i j}$. We define a tensor $\alpha (A) \in {\calT}_4 V$ by
\begin{eqnarray}
\alpha (A)_{i j k l} & := &
{\textstyle \frac{1}{3}} \left( 2\,A_{i j} A_{k l} + A_{i k} A_{j l}
-  A_{i l} A_{j k} \right) \,. \label{alphaa}
\end{eqnarray}
}
\end{enumerate}
\end{Def}
Tensors of type $\gamma (S)$ and $\alpha (A)$ occur in formulas of so-called {\itshape Osserman tensors} (see \cite{gilkey3}).
\begin{Def}
An algebraic curvature tensor $T$ is called {\itshape spacelike Osserman} (resp. {\itshape timelike Osserman}) if the eigenvalues of $J_T(x)$ are constant on $S^{+}(V):=\{x\in V\,|\,g(x,x)=+1\}$ (resp. $S^{-}(V):=\{x\in V\,|\,g(x,x)=-1\}$).
\end{Def}
Since the fundamental tensor $g$ connects an algebraic curvature tensor $T$ with its Jacobi operator $J_T(x)$ in (\ref{jacobi}), the Osserman property depends not only on $T$ but also on $g$. It is known that spacelike Osserman and timelike Osserman are equivalent notions so one simply says {\itshape Osserman}.

If $R$ is the Riemann tensor of a Riemannian manifold $(M,g)$ which is locally a rank one symmetric space or flat, then the eigenvalues of $J_R(x)$ are constant on the unit sphere. Osserman \cite{oss90} wondered if the converse held. This question is known as the {\itshape Osserman conjecture}.

%First page headline in AmS-LaTeX for S\'eminaire Lotharingien de Combinatoire
%--restoring the headers and pagenumbering
%\pagenumbering{arabic}
%\addtocounter{page}{1}
%\markboth{\SMALL B. FIEDLER}{\SMALL ALGEBRAIC CURVATURE TENSORS}
%
%

The correctness of the Osserman conjecture has been established for Riemannian manifolds $(M,g)$ in all dimensions $\not= 8, 16$ (see \cite{c88,n02}) and for Lorentzian manifolds $(M,g)$ in all dimensions (see \cite{bbg97,gkv97}). The situation is much more complicated in the case of a pseudo-Riemannian metric with signature $(p,q)$, $p,q\ge 2$. It is known that there exist pseudo-Riemannian Osserman manifolds which are not locally homogeneous \cite{bbgr97,bcghv98,gvv98}. Furthermore there are counter examples to the conjecture if the manifold $(M,g)$ has signature $(2,2)$ and if the Jacobi operator is not diagonalizable \cite{bbgr97,bbr98}. (See also the examples in \cite{gkvv99,ra97}). Thus the question Osserman raised has a negative answer in the higher signature setting. 
More extensive bibliographies in the subject can be found in the books \cite{gilkey5} by P. Gilkey and \cite{gkv02} by E. Garc\'ia-Rio, D. N. Kupeli and R. V\'azquez-Lorenzo.

Examples of Osserman algebraic curvature tensors can be constructed in the following way. Let $\{ C_i \}_{i = 1}^r$ be a finite set of real, skew-symmetric
$(\mathrm{dim} V \times \mathrm{dim} V)$-matrices that satisfy the Clifford commutation relations
\begin{eqnarray}
C_i \cdot C_j + C_j \cdot C_i & = & - 2\,{\delta}_{i j} \,. \label{equ1.8}
\end{eqnarray}
We define
\begin{eqnarray}
T_C (x, y) z & := &
g(y, C\cdot z) \,C\cdot x -
g(x, C\cdot z) \,C\cdot y -
2\,g(x, C\cdot y) \,C\cdot z \,, \label{equ1.9}\\
T_0 (x, y) z & := &
g(y, z) \,x -
g(x, z) \,y \,, \;\;\;\;\; x, y, z \in V\,,\\
T & := &
{\lambda}_0 \,T_0 +
{\lambda}_1 \,T_{C_1} + \ldots +
{\lambda}_r \,T_{C_r} \,. \label{osserman1}
\end{eqnarray}
Now if $g$ is positive definite, then $T$ is an algebraic curvature tensor that is {\itshape Osserman}. (See \cite{gilkey1,gls99,gilkey3} for more details.)

$T_0$ and $T_C$ can be expressed by means of $\gamma$ and $\alpha$, respectively, since
\begin{eqnarray*}
g(T_0(x, y)\,z, w) & = &
g(y, z) \,g(x, w) - g(x, z) \,g(y, w) \\
 & = & \left( g_{j k} g_{i l} - g_{i k} g_{j l} \right) x^i y^j z^k w^l \\
 & = & 3\, {\gamma}(g)_{i j k l}\,  x^i y^j z^k w^l
\end{eqnarray*}
and
\begin{eqnarray*}
g(T_C(x, y)\,z, w) & = &
g(y, C\cdot z) \,g(C\cdot x, w) -
g(x, C\cdot z) \,g(C\cdot y, w) -\\
 & & 2\,g(x, C\cdot y) \,g(C\cdot z, w)\\
 & = &
g_{j s} y^j {C^s}_k z^k\,g_{t l} {C^t}_i x^i w^l -
g_{i s} x^i {C^s}_k z^k\,g_{t l} {C^t}_j y^j w^l -\\
 & & 
2\,g_{i s} x^i {C^s}_j y^j\,g_{t l} {C^t}_k z^k w^l \\
 & = &
\left\{
C_{j k} C_{l i} -
C_{i k} C_{l j} - 2\,C_{i j} C_{l k} \right\}\,
x^i y^j z^k w^l\\
 & = &
\left\{
- C_{j k} C_{i l} +
C_{i k} C_{jl } + 2\,C_{i j} C_{kl } \right\}\,
x^i y^j z^k w^l\\
 & = & 3\, {\alpha}(C)_{i j k l}\,x^i y^j z^k w^l \,.
\end{eqnarray*}
Thus we obtain
\begin{eqnarray}
T & = & 3 \,{\lambda}_0\,\gamma (g) + 3 \sum_{i = 1}^r
{\lambda}_i \,\alpha (C_i) \label{osserman2}
\end{eqnarray}
from \eqref{osserman1}. It was conjectured that any algebraic Osserman curvature tensor has a structure given in \eqref{osserman1}. This is true for all known examples of Osserman tensors in the case of a positive definite metric $g$ (see \cite{gilkey3}). But for indefinite metrics there exist counter examples (see Remark \ref{rem6.7}).

In the present paper we show that every algebraic curvature tensor $T$ has a representation
\begin{eqnarray}
T & = & \sum_{S \in \calS} {\epsilon}_S\,\gamma (S) \;+\; \sum_{A \in \calA} {\epsilon}_A\,\alpha (A)
\;\;\;\;,\;\;\;\;
{\epsilon}_S\,,\,{\epsilon}_A \in \{1 \;,\; -1\}\,, \label{algcurv}
\end{eqnarray}
where $\calS$ and $\calA$ are finite sets of symmetric or skew-symmetric tensors of order 2, respectively. (See Theorem \ref{algcurvthm}.) To prove \eqref{algcurv} we use the connection between tensors and the representation theory of symmetric groups. We discussed the formula \eqref{algcurv} with P. B. Gilkey. In the course of this discussion P. B. Gilkey found another proof of \eqref{algcurv} which is based on tensor algebra (see \cite[pp.41-44]{gilkey5}). More precisely, P. B. Gilkey could show the stronger result that every algebraic curvature tensor has a representation
\begin{eqnarray}
T & = & \sum_{S \in \calS} {\epsilon}_S\,\gamma (S)
\;\;\;\;,\;\;\;\;
{\epsilon}_S \in \{1 \;,\; -1\}\,, \label{algcurvs}
\end{eqnarray}
as well as a representation
\begin{eqnarray}
T & = & \sum_{A \in \calA} {\epsilon}_A\,\alpha (A)
\;\;\;\;,\;\;\;\;
{\epsilon}_A \in \{1 \;,\; -1\}\,. \label{algcurva}
\end{eqnarray}
In Section \ref{sec5} we give a new proof of \eqref{algcurvs} and \eqref{algcurva} using symmetry operators.

In Section \ref{sec4} we verify our results by means of the Littlewood-Richardson rule and plethysms. We finish the paper with some remarks about applications of tensors $\gamma(S)$ and $\alpha(A)$ in the theory of Osserman algebraic curvature tensors.
\vspace{1cm}
\section{The symmetry class of algebraic curvature tensors}
The set of the algebraic curvature tensors over $V$ is a symmetry class in the sence of H. Boerner \cite[p.127]{boerner}. We denote by $\bbK [{\calS}_r]$ the {\itshape group ring} of a symmetric group ${\calS}_r$. Every group ring element $a = \sum_{p \in {\calS}_r} a(p)\,p \in \bbK [{\calS}_r]$ acts as so-called {\itshape symmetry operator} on tensors $T \in {\calT}_r V$ according to the definition
\begin{eqnarray}
(a T)(v_1 , \ldots , v_r) & := & \sum_{p \in {\calS}_r} a(p)\,
T(v_{p(1)}, \ldots , v_{p(r)}) \;\;\;\;\;,\;\;\;\;\;
v_i \in V \,. \label{symop}
\end{eqnarray}
Equation \eqref{symop} is equivalent to
\begin{eqnarray}
(a T)_{i_1 \ldots i_r} & = & \sum_{p \in {\calS}_r} a(p)\,
T_{i_{p(1)} \ldots  i_{p(r)}} \,.
\end{eqnarray}
\begin{Def}
Let $\frakr \subseteq \bbK [{\calS}_r]$ be a right ideal of $\bbK [{\calS}_r]$ for which an $a \in \frakr$ and a $T \in {\calT}_r V$ exist such that $aT \not= 0$. Then the tensor set
\begin{eqnarray}
{\calT}_{\frakr} & := & \{ a T \;|\; a \in \frakr \;,\;
T \in {\calT}_r V \}
\end{eqnarray}
is called the {\itshape symmetry class} of tensors defined by $\frakr$.
\end{Def}
Since $\bbK [{\calS}_r]$ is semisimple for $\bbK = \bbR , \bbC$, every right ideal $\frakr \subseteq \bbK [{\calS}_r]$ possesses a generating idempotent $e$, i.e. $\frakr$ fulfils $\frakr = e \cdot \bbK [{\calS}_r]$.
\begin{Lem}
If $e$ is a generating idempotent of $\frakr$, then a tensor $T \in {\calT}_r V$ belongs to ${\calT}_{\frakr}$ iff
\begin{eqnarray}
e T & = & T \,.
\end{eqnarray}
Thus we have
\begin{eqnarray}
{\calT}_{\frakr} & = & \{ eT \;|\; T \in {\calT}_r V \} \,. \label{classgen}
\end{eqnarray}
\end{Lem}
\begin{proof}
Since every ideal element $a \in \frakr$ can be written as
$a = e \cdot x$ with $x \in \bbK [{\calS}_r]$, every $aT \in {\calT}_r V$ fulfils $aT = e(xT) = e T'$ with $T' = xT \in {\calT}_r V$. Thus we obtain
${\calT}_r \subseteq \{ eT \;|\; T \in {\calT}_r V \}$. The relation
$\{ eT \;|\; T \in {\calT}_r V \} \subseteq {\calT}_r$ is trivial.
\end{proof}
A {\itshape Young frame} of $r$ is an arrangement of $r$ boxes such that the numbers ${\lambda}_i$ of boxes in the rows $i = 1 , \ldots , l$ form a decreasing sequence
${\lambda}_1 \ge {\lambda}_2 \ge \ldots \ge {\lambda}_l > 0$ with
${\lambda}_1 + \ldots + {\lambda}_l = r$.
\[
\begin{array}{cc|c|c|c|c|c|c}
\cline{3-7}
{\lambda}_1 = 5 & \;\;\; & \;\;\; & \;\;\; & \;\;\; & \;\;\; & \;\;\; & \\
\cline{3-7}
{\lambda}_2 = 4 & \;\;\; & \;\;\; & \;\;\; & \;\;\; & \;\;\; & \multicolumn{2}{c}{\;\;\;} \\
\cline{3-6}
{\lambda}_3 = 4 & \;\;\; & \;\;\; & \;\;\; & \;\;\; & \;\;\; & \multicolumn{2}{c}{\;\;\;} \\
\cline{3-6}
{\lambda}_4 = 2 & \;\;\; & \;\;\; & \;\;\; & \multicolumn{4}{c}{\hspace{2cm}} \\
\cline{3-4}
{\lambda}_5 = 1 & \;\;\; & \;\;\; & \multicolumn{4}{c}{\hspace{2cm}} \\
\cline{3-3}
\end{array}
\]
Obviously, a Young frame of $r$ is characterized by a partition
$\lambda = ({\lambda}_1 , \ldots , {\lambda}_l) \vdash r$ of $r$. A {\itshape Young tableau} $t$ of a partition $\lambda \vdash r$ is the Young frame corresponding to $\lambda$ which was fulfilled by the numbers $1, 2, \ldots , r$ in any order.
\begin{eqnarray*}
t & = &
\begin{array}{|c|c|c|c|c|c}
\cline{1-5}
11 & 2 & 5 & 4 & 12 & \\
\cline{1-5}
9 & 6 & 16 & 15 & \multicolumn{2}{c}{\;\;\;} \\
\cline{1-4}
8 & 14 & 1 & 7 & \multicolumn{2}{c}{\;\;\;} \\
\cline{1-4}
13 & 3 & \multicolumn{4}{c}{\hspace{2cm}} \\
\cline{1-2}
10 & \multicolumn{4}{c}{\hspace{2cm}} \\
\cline{1-1}
\end{array}
\end{eqnarray*}
If a Young tableau $t$ of a partition $\lambda \vdash r$ is given, then the {\itshape Young symmetrizer} $y_t$ of $t$ is defined by\footnote{We use the convention $(p \circ q) (i) := p(q(i))$ for the product of two permutations $p, q$.}
\begin{eqnarray}
y_t & := & \sum_{p \in {\calH}_t} \sum_{q \in {\calV}_t} \mathrm{sign}(q)\, p \circ q
\end{eqnarray}
where ${\calH}_t$, ${\calV}_t$ are the groups of the {\itshape horizontal} or
{\itshape vertical permutations} of $t$ which only permute numbers within rows or columns of $t$, respectively. The Young symmetrizers of $\bbK [{\calS}_r]$ are essentially idempotent and define decompositions
\begin{eqnarray}
\bbK [{\calS}_r] \;=\;
\bigoplus_{\lambda \vdash r} \bigoplus_{t \in {\calS\calT}_{\lambda}}
\bbK [{\calS}_r]\cdot y_t
& \;\;,\;\; &
\bbK [{\calS}_r] \;=\;
\bigoplus_{\lambda \vdash r} \bigoplus_{t \in {\calS\calT}_{\lambda}}
y_t \cdot \bbK [{\calS}_r] \label{decomp}
\end{eqnarray}
of $\bbK [{\calS}_r]$ into minimal left or right ideals. In \eqref{decomp}, the symbol ${\calS\calT}_{\lambda}$ denotes the set of all standard tableaux of the partition $\lambda$. Standard tableaux are Young tableaux in which the entries of every row and every column form an increasing number sequence.\footnote{About Young symmetrizers and
Young tableaux see for instance
\cite{boerner,boerner2,full4,fulton,jameskerb,kerber,littlew1,mcdonald,%
muell,naimark,%
waerden,weyl1}. In particular, properties of Young symmetrizers in the case
${\bbK} \not= {\bbC}$ are described in \cite{muell}.}

S.A. Fulling, R.C. King, B.G.Wybourne and C.J. Cummins showed in \cite{full4} that the symmetry classes of the Riemannian curvature tensor $R$ and its covariant derivatives ${\nabla}^{(u)}R$ are generated by special  Young symmetrizers. To this end they assumed that the covariant derivatives
${\nabla}_{s_1} {\nabla}_{s_2} \ldots {\nabla}_{s_u} R_{i j k l}$ are symmetric in $s_1, \ldots , s_u$ for $u \ge 2$. This is possible if one uses Ricci identities and neglect derivatives of orders smaller than $u$. Under this assumption one has
\begin{Thm}
Consider the Levi-Civita connection $\nabla$ of a pseudo-Riemannian metric $g$.
For $u \ge 0$ the Riemann tensor and its covariant derivatives
${\nabla}^{(u)} R$ fulfil
\begin{eqnarray}
e_t^{\ast} {\nabla}^{(u)} R & = & {\nabla}^{(u)} R
\end{eqnarray}
where $e_t := y_t (u+1)/(2\cdot (u+3)!)$ is an idempotent which is formed from the Young symmetrizer $y_t$ of the standard tableau
\begin{eqnarray}
t & = &
\begin{array}{|c|c|c|cc|c|}
\hline
1 & 3 & 5 & \ldots & \ldots & (u+4) \\
\hline
2 & 4 & \multicolumn{4}{l}{\hspace{3cm}} \\
\cline{1-2}
\end{array} \,.
\end{eqnarray}
The star ''$\ast$'' denotes the mapping
\begin{eqnarray}
\ast : a \;=\; \sum_{p \in {\calS}_r} a(p)\,p & \mapsto &
a^{\ast} \;:=\; \sum_{p \in {\calS}_r} a(p)\,p^{-1} \,.
\end{eqnarray} \label{riemclass}
\end{Thm}
A proof of this result of \cite{full4} can be found in \cite{fie12}, too. The proof needs only the symmetry properties \eqref{riemsym}, the identities Bianchi I and Bianchi II and Ricci identities. Thus Theorem \ref{riemclass} is a statement about the algebraic curvature tensors and ''algebraic'' tensors that possess the symmetry properties of the covariant derivatives of the Riemann tensor. In particular, a tensor $T \in {\calT}_4 V$ is an algebraic curvature tensor iff $T$ satisfies
\begin{eqnarray}
y_t^{\ast} T & = & 12\,T \label{acgen}
\end{eqnarray}
where $y_t$ is the Young symmetrizer of the standard tableau
\begin{eqnarray}
t & = &
\begin{array}{|c|c|}
\hline
1 & 3 \\
\hline
2 & 4 \\
\hline
\end{array} \,. \label{riemtabl}
\end{eqnarray}
\vspace{1cm}
\section{A structure theorem for algebraic curvature tensors}
Now we will show that every algebraic curvature tensor can be built from certain finite sets of symmetric and skew-symmetric tensors from ${\calT}_2 V$.
\begin{Thm} \label{algcurvthm}
For every algebraic curvature tensor $T \in {\calT}_4 V$ there exist finite sets $\calS$ and $\calA$ of symmetric and alternating tensors $S , A \in {\calT}_2 V$, respectively, such that
\begin{eqnarray}
T & = & \sum_{S \in \calS} {\epsilon}_S\,\gamma (S) \;+\; \sum_{A \in \calA} {\epsilon}_A\,\alpha (A)
\;\;\;\;,\;\;\;\;
{\epsilon}_S\,,\,{\epsilon}_A \in \{1 \;,\; -1\}\,. \label{algcurvrep}
\end{eqnarray}
\end{Thm}
\begin{proof}
Let $T$ be an algebraic curvature tensor.
Because of the symmetry $T_{k l i j} = T_{i j k l}$, every algebraic curvature tensor $T$ fulfils $fT = T$, where $f$ is the idempotent
$f := \frac{1}{2} ( \mathrm{id} + (1\,3)(2\,4)) \in \bbK [{\calS}_4]$. Thus every algebraic curvature tensor belongs to the symmetry class ${\calT}_{\frakr'}$ of the right ideal $\frakr' := f \cdot \bbK [{\calS}_4]$ and has a representation $T = f T'$ with $T' \in {\calT}_4 V$.

Since an arbitrary tensor $T' \in {\calT}_4 V$ can be written as a finite sum of decomposable tensors $a_1 \otimes a_2 \otimes a_3 \otimes a_4 \in {\calT}_4 V$, we can use a representation
\begin{eqnarray}
T' & = & \sum_{(M, N) \in \calP} M \otimes N \label{matrixsum}
\end{eqnarray}
for a $T' \in {\calT}_4 V$, where $\calP$ is a finite set of pairs $(M,N)$ of tensors $M, N \in {\calT}_2 V$.
From \eqref{matrixsum} and $T = fT'$ we obtain
\begin{eqnarray*}
T \;=\; fT' & = & \frac{1}{2}\,\sum_{(M,N) \in \calP}
(M \otimes N + N \otimes M) \\
 & = &
\frac{1}{2}\,\sum_{(M,N) \in \calP}
\left( (M + N) \otimes (M + N) - M \otimes M - N \otimes N \right) \,.
\end{eqnarray*}
Thus every tensor $T \in {\calT}_{\frakr'}$ has a representation
\begin{eqnarray}
T & = & \frac{1}{2}\,\sum_{M' \in {\calM}'} z_{M'}\,M' \otimes M'
\;\;\;\;,\;\;\;\;
z_{M'} \in \bbZ
\label{mmsum1}
\end{eqnarray}
with a certain finite set ${\calM}'$ of tensors from ${\calT}_2 V$.
If we carry out a transformation $M' \mapsto M := \sqrt{|z_{M'}|/2}\, M'$ in every summand of \eqref{mmsum1}, we obtain
\begin{eqnarray}
T & = & \sum_{M \in \calM} {\epsilon}_M\,M \otimes M
\;\;\;\;,\;\;\;\;
{\epsilon}_M \in \{1 \;,\; -1\}
\label{mmsum}
\end{eqnarray}
from \eqref{mmsum1}.

Since every matrix $M$ has a unique decomposition into a symmetric matrix $S$ and a skew-symmetric matrix $A$, i.e. $M = S + A$, \eqref{mmsum} leads to
\begin{eqnarray}
T & = & 
\sum_{S+A \in \calM}
{\epsilon}_{S+A}\,\left( S \otimes S + S \otimes A + A \otimes S + A \otimes A \right) \,. \label{sasum}
\end{eqnarray}

Now, let $y_t$ be the Young symmetrizer of the standard tableau \eqref{riemtabl}. If we apply \eqref{acgen} to the representation \eqref{sasum} of an algebraic curvature tensor $T$, we obtain
\begin{eqnarray}
\;\;\;\;\;T \;=\; {\textstyle \frac{1}{12}}\,y_t^{\ast} T & = &
{\textstyle \frac{1}{12}}\,\sum_{S+A \in \calM}
{\epsilon}_{S+A}\;y_t^{\ast} \left( S \otimes S + S \otimes A + A \otimes S + A \otimes A \right) \,.
\end{eqnarray}
Then a calculation\footnote{This calculation can be done by hand or by a cooperation of the \textsf{Mathematica} package \textsf{PERMS} \cite{fie10} with one of the tensor packages \textsf{Ricci} \cite{ricci3} or \textsf{MathTensor} \cite{mathtensor3}. See \cite{fie21}.} yields
\begin{eqnarray*}
{\textstyle \frac{1}{12}}\,y_t^{\ast} (S \otimes S) \;=\; \gamma (S)
& \;\;\;\;\;\;,\;\;\;\;\;\; & 
{\textstyle \frac{1}{12}}\,y_t^{\ast} (A \otimes S) \;=\; 0\,, \\
{\textstyle \frac{1}{12}}\,y_t^{\ast} (A \otimes A) \;=\; \alpha (A)
& \;\;\;\;\;\;,\;\;\;\;\;\; & 
{\textstyle \frac{1}{12}}\,y_t^{\ast} (S \otimes A) \;=\; 0 \,,
\end{eqnarray*}
and if we set ${\epsilon}_S := {\epsilon}_{S+A}$, ${\epsilon}_A := {\epsilon}_{S+A}$, we obtain \eqref{algcurvrep}.
\end{proof}
\begin{Cor}
If $\calS \subseteq {\calT}_2 V$ and $\calA \subseteq {\calT}_2 V$ are finite sets of symmetric or alternating tensors, respectively, then
\begin{eqnarray}
T & := & \sum_{S \in \calS} {\epsilon}_S\,\gamma (S) \;+\; \sum_{A \in \calA} {\epsilon}_A\,\alpha (A)
\;\;\;\;,\;\;\;\;
{\epsilon}_S\,,\,{\epsilon}_A \in \{1 \;,\; -1\} \label{algcurvrep2}
\end{eqnarray}
is an algebraic curvature tensor.
\end{Cor}
\begin{proof}
Every summand of \eqref{algcurvrep2} is an algebraic curvature tensor since
$\gamma (S) = \frac{1}{12}\,y_t^{\ast} (S \otimes S)$ and
$\alpha (A) = \frac{1}{12}\,y_t^{\ast} (A \otimes A)$.
\end{proof}
\vspace{1cm}
\section{Verification by the Littlewood-Richardson rule and plethysms} \label{sec4}
We can verify the results about
$y_t^{\ast} (S \otimes S)$,
$y_t^{\ast} (S \otimes A)$,
$y_t^{\ast} (A \otimes S)$ and
$y_t^{\ast} (A \otimes A)$ by means of the \textit{Littlewood-Richardson rule}\footnote{See
D. E. Littlewood \cite[pp.94-96]{littlew1},
A. Kerber \cite[Vol.240/p.84]{kerber},
G. D. James and A. Kerber \cite[p.93]{jameskerb},
A. Kerber \cite[Sec.5.5]{kerber3},
I. G. Macdonald \cite[Chap.I,Sec.9]{mcdonald},
W. Fulton and J. Harris \cite[pp.455-456]{fultharr},
S. A. Fulling, R. C. King, B. G. Wybourne and C. J. Cummins \cite{full4}. See also B. Fiedler \cite[Sec.II.5]{fie16}.} and \textit{plethysms}\footnote{See
A. Kerber \cite{kerber}, G. D. James and A. Kerber 
\cite{jameskerb}, A. Kerber \cite{kerber3}, F. S\"anger \cite{saenger},
D. E. Littlewood
\cite{littlew1}, I. G. Macdonald \cite{mcdonald},
S. A. Fulling, R. C. King, B. G. Wybourne and C.J. Cummins \cite{full4}. See also B. Fiedler \cite[Sec.II.6]{fie16}.}. To this end we form group ring elements from tensors (see B. Fiedler \cite[Sec.III.1]{fie16} and \cite{fie8,fie14,fie17}).
\begin{Def}
Any tensor $T \in {\calT}_r V$ and any $r$-tuple
$b := (v_1 , \ldots , v_r) \in V^r$ of vectors from $V$ induce a function
$T_b : {\calS}_r \rightarrow \bbK$ according to the rule
\begin{eqnarray}
T_b (p) & := & T(v_{p(1)}, \ldots , v_{p(r)}) \;\;\;,\;\;\;
p \in {\calS}_r \,.
\end{eqnarray}
We identify this function with the group ring element
$\sum_{p \in {\calS}_r} T_b(p)\,p \in \bbK [{\calS}_r]$,
which we denote by $T_b$, too.
\end{Def}
These $T_b$ fulfil$^{\ref{foot}}$
\begin{Cor}
If $T \in {\calT}_r V$, $b \in V^r$ and $a \in \bbK [{\calS}_r]$, then
\begin{eqnarray}
(aT)_b & = & T_b \cdot a^{\ast} \,. \label{tbrule}
\end{eqnarray}
\end{Cor}
If a tensor $T$ belongs to a certain symmetry class, then its $T_b$ lie in a certain left ideal\footnote{See B. Fiedler \cite[Sec.III.3.1]{fie16} and \cite{fie8,fie14,fie17}. \label{foot}}.
\begin{Prop}
Let $e \in \bbK [{\calS}_r]$ be an idempotent. Then a
$T \in {\calT}_r V$ fulfils the condition $eT = T$ iff
$T_b \in \frakl := \bbK [{\calS}_r] \cdot e^{\ast}$ for all
$b \in V^r$, i.e. all $T_b$ of $T$ lie in the left ideal $\frakl$ generated by $e^{\ast}$.
\end{Prop}
Obviously, the $S_b$ and $A_b$ of symmetric/skew-symmetric tensors of order 2 lie in the left ideals
\begin{eqnarray}
\forall\,b \in V^2 : \;\; S_b \in {\frakl}_1 := \bbK [{\calS}_2] \cdot e_1
& \;\;\;,\;\;\; &
e_1 \;:=\; {\textstyle \frac{1}{2}}\,\left(\mathrm{id} + (1\,2) \right) \\
\forall\,b \in V^2 : \;\; A_b \in {\frakl}_2 := \bbK [{\calS}_2] \cdot e_2
& \;\;\;,\;\;\; &
e_2 \;:=\; {\textstyle \frac{1}{2}}\,\left(\mathrm{id} - (1\,2) \right)
\end{eqnarray}
$e_1$ and $e_2$ are the (normalized) Young symmetrizers of the standard tableaux
\begin{eqnarray*}
t_1 \;:=\;
\begin{array}{|c|c|}
\hline
1 & 2 \\
\hline
\end{array}
& \;\;\;\;\;,\;\;\;\;\; &
t_2 \;:=\;
\begin{array}{|c|}
\hline
1 \\
\hline
2 \\
\hline
\end{array}
\end{eqnarray*}
and the left ideals ${\frakl}_1$ and ${\frakl}_2$ are 1-dimensional, minimal left ideals, which belong to the equivalence classes of minimal left ideals characterized by the partitions ${\lambda}_1 = (2)$ and ${\lambda}_2 = (1^2)$, respectively. Then the group ring elements
$(A \otimes S)_b$,
$(S \otimes A)_b$,
$(S \otimes S)_b$ and
$(A \otimes A)_b$ belong to left ideals which are characterized by Littlewood-Richardson products $[2] [1^2]$ or plethysms $[2] \odot [2]$, $[1^2] \odot [2]$ (see B. Fiedler \cite[Sec.III.3.2]{fie16} or \cite[Sec.4.2]{fie17}). Especially, we obtain the following structures for these ideals:
\[
\begin{array}{|c|c|l|}
\hline
(S \otimes S)_b &  [2] \odot [2]  \;\sim\;  [4] + [2^2] &
{\frakl}^{(1)} \;=\; {\frakl}_{[4]}^{(1)} \oplus {\frakl}_{[2^2]}^{(1)}\\
\hline
(S \otimes A)_b &  [2][1^2]  \;\sim\;  [3\,1] + [2\,1^2] &
{\frakl}^{(2)} \;=\; {\frakl}_{[3\,1]}^{(2)} \oplus {\frakl}_{[2\,1^2]}^{(2)}\\
\hline
(A \otimes S)_b &  [2][1^2]  \;\sim\;  [3\,1] + [2\,1^2] &
{\frakl}^{(3)} \;=\; {\frakl}_{[3\,1]}^{(3)} \oplus {\frakl}_{[2\,1^2]}^{(3)}\\
\hline
(A \otimes A)_b &  [1^2] \odot [2]  \;\sim\;  [2^2] + [1^4] &
{\frakl}^{(4)} \;=\; {\frakl}_{[2^2]}^{(4)} \oplus {\frakl}_{[1^4]}^{(4)}\\
\hline
\end{array}
\]
The relation $[2][1^2] \sim [3\,1] + [2\,1^2]$ can be determined by the Littlewood-Richardson rule. The relation $[2] \odot [2] \sim  [4] + [2^2]$ follows from the formula
\begin{eqnarray}
[2] \odot [n] & \sim & \sum_{\lambda \vdash n} \,[2\,\lambda] \label{pleth2n}
\end{eqnarray}
(see G.D. James and A. Kerber \cite[p.224]{jameskerb}) and
$[1^2] \odot [2] \sim  [2^2] + [1^4]$ is a consequence of \eqref{pleth2n} and
\begin{eqnarray}
[\mu] \odot [\nu] \;\sim\; \sum_{\lambda \vdash m\cdot n} m_{\lambda} [\lambda']
& \Leftrightarrow &
[\mu'] \odot [\nu] \;\sim\; \sum_{\lambda \vdash m\cdot n} m_{\lambda} [\lambda] \label{transpleth} \\
 & & \mu \vdash m \,,\, \nu \vdash n \,,\,m, n, m_{\lambda} \in \bbN \,,\,
m\,\mathrm{even}\,. \nonumber
\end{eqnarray}
(See F. S\"anger \cite[pp.10-11]{saenger}. See also B. Fiedler \cite[Corr.II.6.8]{fie16}.) In \eqref{transpleth} $\lambda'$ denotes the associated or transposed partition of a partition $\lambda$. All formulas of the above table can be determined by the  \textsf{Mathematica} package \textsf{Perms} \cite{fie10}, too.

Every of the left ideals ${\frakl}^{(k)}$, $k = 1, 2, 3, 4$, decomposes into two minimal left ideals
${\frakl}^{(k)} = {\frakl}_{[{\mu}_k]}^{(k)} \oplus {\frakl}_{[{\nu}_k]}^{(k)}$, where the partitions ${\mu}_k$ and ${\nu}_k$ charakterize the equivalence class of minimal left ideals to which such an ideal 
${\frakl}_{[{\mu}_k]}^{(k)}\,,\,{\frakl}_{[{\nu}_k]}^{(k)}$ belongs.

Now the formula \eqref{tbrule}
yields that the group ring elements
$(y_t^{\ast}(S \otimes S))_b$,
$(y_t^{\ast}(A \otimes S))_b$,
$(y_t^{\ast}(S \otimes A))_b$ and
$(y_t^{\ast}(A \otimes A))_b$
are contained in one of the left ideals ${\frakl}^{(k)} \cdot y_t$. Since the Young symmetrizer $y_t$ of \eqref{riemtabl} generates a minimal left ideal of the equivalence class of $\lambda = (2^2) \vdash 4$, we obtain
${\frakl}_{[\mu]}^{(k)} \cdot y_t = \{ 0 \}$ if $\mu \not= (2^2)$. This leads to ${\frakl}^{(2)} \cdot y_t = \{ 0 \}$, ${\frakl}^{(3)} \cdot y_t = \{ 0 \}$ and
$y_t^{\ast}(A \otimes S) = 0$,
$y_t^{\ast}(S \otimes A)$. Furthermore, we obtain
${\frakl}_{[4]}^{(1)} \cdot y_t = \{ 0 \}$,
${\frakl}_{[1^4]}^{(4)} \cdot y_t = \{ 0 \}$. At most the products
${\frakl}_{[2^2]}^{(1)} \cdot y_t$,
${\frakl}_{[2^2]}^{(4)} \cdot y_t$ are candidates for non-vanishing results. We have already seen in the proof of Theorem \ref{algcurvthm} that these cases really lead to non-vanishing expressions \eqref{gammas} and \eqref{alphaa}.
\vspace{1cm}
\section{A new proof of the stronger result of P. B. Gilkey} \label{sec5}
\noindent Now we give a new proof of the following result of P. B. Gilkey (see \cite[pp.41-44]{gilkey5}).
\begin{Thm} \label{gilkeythm}
For every algebraic curvature tensor $T \in {\calT}_4 V$ there exist finite sets $\calS$ and $\calA$ of symmetric and alternating tensors $S , A \in {\calT}_2 V$, respectively, such that $T$ has a representation
\begin{eqnarray}
T & = & \sum_{S \in \calS} {\epsilon}_S\,\gamma (S)
\;\;\;\;,\;\;\;\;
{\epsilon}_S \in \{1 \;,\; -1\}\,, \label{algcurvreps}
\end{eqnarray}
as well as a representation
\begin{eqnarray}
T & = & \sum_{A \in \calA} {\epsilon}_A\,\alpha (A)
\;\;\;\;,\;\;\;\;
{\epsilon}_A \in \{1 \;,\; -1\}\,. \label{algcurvrepa}
\end{eqnarray}
\end{Thm}
\begin{proof}
We consider the following group ring elements
\begin{eqnarray}
{\sigma}_{\epsilon} & := & y_t^{\ast} \cdot f \cdot {\xi}_{\epsilon}
\;\;\;\;,\;\;\;\;
{\epsilon} \in \{1 \;,\; -1\}\,, \label{sym}
\end{eqnarray}
where $y_t$ is the Young symmetrizer of the Young tableau \eqref{riemtabl} and
\begin{eqnarray}
f & := & \mathrm{id} + (1\,3)(2\,4) \label{symf}\\
{\xi}_{\epsilon} & := &
\bigl( \mathrm{id} + \epsilon\,(1\,2)\bigr) \cdot 
\bigl( \mathrm{id} + \epsilon\,(3\,4)\bigr)
\;\;\;\;,\;\;\;\;
{\epsilon} \in \{1 \;,\; -1\}\,. \label{symxi}
\end{eqnarray}
A calculation (by means of {\sf PERMS} \cite{fie10}) shows that ${\sigma}_1 \not= 0$ and ${\sigma}_{-1} \not= 0$ (see \cite{fie21}).

The right ideals ${\sigma}_{\epsilon} \cdot \bbK [\calS_4]$ are non-vanishing subideals of the right ideal $\frakr = y_t^{\ast}\cdot \bbK [\calS_4]$ which defines the symmetry class ${\calT}_{\frakr}$ of algebraic curvature tensors. Since $\frakr$ is a minimal right ideal, we obtain
$\frakr = {\sigma}_1 \cdot \bbK [\calS_4] = {\sigma}_{-1} \cdot \bbK [\calS_4]$, i.e. ${\sigma}_1$ and ${\sigma}_{-1}$ are generating elements of $\frakr$, too.

A tensor $T \in {\calT}_4 V$ is an algebraic curvature tensor iff there exist $a \in \frakr$ and $T' \in {\calT}_4 V$ such that $T = aT'$. Since further every $a \in \frakr$ has a representation
$a = {\sigma}_{\epsilon}\cdot x_{\epsilon}$, $\epsilon = 1 , -1$, a tensor
$T \in {\calT}_4 V$ is an algebraic curvature tensor iff there is a tensor
$T' \in {\calT}_4 V$ such that $T = {\sigma}_1 T'$ or
$T = {\sigma}_{-1} T'$.

Let us consider the case $\epsilon = 1$. We obtain all algebraic curvature tensors if we form $T = {\sigma}_1 T'$, $T' \in {\calT}_4 V$. As in the proof of Theorem \ref{algcurvthm} we can use a representation
\begin{eqnarray}
T' & = & \sum_{(M,N) \in \calP} M \otimes N \label{matrixsum2}
\end{eqnarray}
for a $T' \in {\calT}_4 V$, where $\calP$ is a finite set of pairs $(M,N)$ of tensors $M , N \in {\calT}_2 V$. From \eqref{matrixsum2} we obtain
\begin{eqnarray}
{\xi}_1 T' & = & \sum_{(S',S'') \in {\tilde \calP}} S' \otimes S'' \label{sum3}
\end{eqnarray}
where ${\tilde \calP}$ is the finite set of pairs $(S',S'')$ of the symmetrized tensors $S' = M + M^t$, $S'' = N + N^t$, $(M,N) \in \calP$.
Further the application of the symmetry operator $f$ to ${\xi}_1 T'$ yields
\begin{eqnarray}
f({\xi}_1 T') & = & \sum_{(S',S'') \in {\tilde \calP}} (S' \otimes S''\,+\, S'' \otimes S' ) \nonumber \\
 & = &
\sum_{(S',S'') \in {\tilde \calP}} \bigl( (S'+S'') \otimes (S'+S'')\,-\, S' \otimes S' \,-\, S'' \otimes S''\bigr) \nonumber \\
 & = &
\sum_{S''' \in \calS'''} z_{S'''}\,S''' \otimes S'''
\;\;\;\;,\;\;\;\;
z_{S'''} \in \bbZ \nonumber \\
 & = &
\sum_{S \in \calS} {\epsilon}_S\,S \otimes S
\;\;\;\;,\;\;\;\;
{\epsilon}_S \in \{1 \;,\; -1\}\,.
\end{eqnarray}
The last step includes a transformation
$S''' \mapsto S := \sqrt{|z_{S'''}|} S'''$. The application of $y_t^{\ast}$ finishs the proof of \eqref{algcurvreps}:
\[
T \;=\;
{\sigma}_1 T'
\;=\;
(y_t^{\ast}\cdot f\cdot {\xi}_1) T'
\;=\;
\sum_{S \in \calS} {\epsilon}_S\,y_t^{\ast}(S \otimes S)
\;=\;
12\,\sum_{S \in \calS} {\epsilon}_S\,\gamma (S) \,.
\]
The proof of \eqref{algcurvrepa} can be carried out in exactly the same way. We have only to replace \eqref{sum3} by
\begin{eqnarray}
{\xi}_{-1} T' & = & \sum_{(A',A'') \in {\hat \calP}} A' \otimes A''
\end{eqnarray}
where ${\hat \calP}$ is the finite set of pairs $(A',A'')$ of the skew-symmetrized tensors $A' = M - M^t$, $A'' = N - N^t$, $(M,N) \in \calP$.
\end{proof}
\begin{Rem}
The new generating elements ${\sigma}_1$, ${\sigma}_{-1}$ of the right ideal $\frakr$ do not automatically lead to new generating idempotents of $\frakr$. A calculation \cite{fie21} by means of {\sf PERMS} \cite{fie10} yields
\begin{eqnarray*}
y_t^{\ast}\cdot y_t^{\ast} & = & 12\, y_t^{\ast}\\
y_t^{\ast}\cdot {\sigma}_1 & = & 12\,{\sigma}_1 \\
{\sigma}_1 \cdot y_t^{\ast} & = & 0 \\
y_t^{\ast}\cdot {\sigma}_{-1} & = & 12\,{\sigma}_{-1} \\
{\sigma}_{-1} \cdot y_t^{\ast} & = & 96\, y_t^{\ast} \\
{\sigma}_1 \cdot {\sigma}_1 & = & 0 \\
{\sigma}_{-1} \cdot {\sigma}_{-1} & = & 96\, {\sigma}_{-1}\\
{\sigma}_{-1} \cdot {\sigma}_{1} & = & 96\, {\sigma}_{1}\\
{\sigma}_{1} \cdot {\sigma}_{-1} & = & 0 \,.
\end{eqnarray*}
Thus ${\sigma}_{-1}$ is an essencially idempotent element from which a generating idempotent of $\frakr$ can be determined. However, ${\sigma}_1$ is a nilpotent generating element of $\frakr$.
\end{Rem}
\vspace{1cm}
\section{Some applications of $\alpha$ and $\gamma$} \label{sec6}
\noindent We finish our paper with some remarks which arise from suggestions and hints of P. B. Gilkey.

In \cite[pp.191-193]{gilkey5} P. B. Gilkey presents a family of algebraic curvature tensors which is a generalization of the family (\ref{osserman1}) in the case of a metric $g$ with an arbitrary signature $(p,q)$. This new family can be used to show that there are {\itshape Jordan Osserman algebraic curvature tensors} $T$ for which $J_T$ can have an arbitrary complicated Jordan normal form. We give a short view of the construction of Gilkey.

We consider a metric $g\in\calT_2 V$ of signature $(p,q)$ on $V$.
\begin{Def}
An algebraic curvature tensor $T\in\calT_4 V$ is called {\itshape spacelike Jordan Osserman} (resp. {\itshape timelike Jordan Osserman}) if the Jordan normal form of the Jacobi operator $J_T(x)$ is constant on $S^{+}(V):=\{x\in V\,|\,g(x,x)=+1\}$ (resp. $S^{-}(V):=\{x\in V\,|\,g(x,x)=-1\}$). $T$ is called {\itshape Jordan Osserman} if $T$ is both spacelike Jordan Osserman and timelike Jordan Osserman.
\end{Def}
Now we generalize the Clifford commutation relation (\ref{equ1.8})
\begin{Def}
Let $V$ and $W$ be vector spaces which possess metrics with signature $(p,q)$ and $(r,s)$, respectively, where $r\le p$ and $s\le q$. We say that $V$ {\itshape admits} a {\itshape Clifford module structure} Cliff($W$) if there exist skew-symmetric linear transformations $C_i : V\rightarrow V$ for $1\le i\le r+s$ so that
\begin{enumerate}
\item{$C_i C_j + C_j C_i = 0$ for $i\not= j$.}
\item{$C_i{}^2 = \Id$ for $r$ values of $i$.}
\item{$C_i{}^2 = - \Id$ for $s$ values of $i$.}
\end{enumerate}
\end{Def}
Note that a linear map $C:V\rightarrow V$ is called {\itshape skew-symmetric} if
\begin{eqnarray*}
\forall\,x,y\in V:\;\;\;g(Cx,y) = - g(x,Cy).
\end{eqnarray*}
Now the following Lemma holds
\begin{Lem} \label{lem6.3}
Let $V$ and $W$ possess metrics of signature $(p,q)$ and $(r,s)$, respectively $(r\le p,\; s\le q)$. Assume that $V$ admits a {\rm Cliff}$(W)$ module structure. Let $J:W\rightarrow W$ be a self-adjoint linear map of $W$. Then there exists a Jordan Osserman algebraic curvature tensor $T$ on $V$ so that the Jacobi operator $J_T(x)$ is conjugate to $g(x,x)\,J\oplus 0$ for every $x\in S^{\pm}(V)$, where $g$ is the metric of $V$.
\end{Lem}
In Gilkey's proof of Lemma \ref{lem6.3} the above Jordan Osserman algebraic curvature tensor $T$ is given by the formula
\begin{eqnarray}
T & = & \sum_i c_i\,\alpha(A_i) + \frac{1}{2}\,\sum_{i\not= j} c_{ij}\,\alpha(A_i + A_j) \label{equ6.1}
\end{eqnarray}
where $c_i$ and $c_{ij} = c_{ji}$ are constants. The $A_i$ are tensors of order 2 on $V$ defined by $A_i(x,y) := g(C_i x,y)$. Obviously the $A_i$ are skew-symmetric.

The expression $J\oplus 0$ is defined as follows: Under the assumtions of Lemma \ref{lem6.3} we can find an isometry between $V$ and $W\oplus\bbR^{(p-r,q-s)}$. Then $J\oplus 0$ is defined by setting
\begin{eqnarray*}
(J\oplus 0)(w\oplus z) & := & Jw\oplus 0.
\end{eqnarray*}
The sum (\ref{equ6.1}) is a representation (\ref{algcurvrepa}) of an algebraic curvature tensor which is formed by means of $\alpha$ from skew-symmetric tensors $A_i\in\calT_2 V$.

Now the following Theorem by Gilkey tells us that we can construct Jordan Osserman algebraic curvature tensors with arbitrarily complicated Jordan normal form if we greatly increase the dimension by suitable powers of 2.
\begin{Thm}
Let $J:W\rightarrow W$ be an arbitrary linear map of a vector space $W$ of dimension $m$. Then there exist $l = l(m)$ and a Jordan Osserman algebraic curvature tensor $T$ on $V = \bbR^{(2^l,2^l)}$ so that $J_T(x)$ is conjugate to $\pm J\oplus 0$ if $x\in S^{\pm}(V)$.
\end{Thm}
If we choose a suitable $J$ we can generate a Jacobi operator $J_T(x)$ with an arbitrarily complicated Jordan normal form. In particular we can construct nondiagonalizable Jacobi operators. Thus the family (\ref{equ6.1}) is different from the family (\ref{osserman1}) since the Jacobi operators of the tensors (\ref{osserman1}) are defined on a vector space $V$ with a positive definite metric and can be diagonalized as self-adjoint operators on $V$.

In \cite[p.191]{gilkey5} P. B. Gilkey presents also a second generalization of (\ref{osserman1}), (\ref{osserman2}) in the case of certain indefinite metrics. The algebraic curvature tensors of this generalization have a structure 
\begin{eqnarray}
T & = & {\lambda}_0\,\gamma(g) + \sum_i {\lambda}_i\,\alpha(A_i)\;\;\;,\;\;\;{\lambda}_0, {\lambda}_i = \mathrm{const.} \label{equ6.2}
\end{eqnarray}
which is very similar to (\ref{osserman2}). The $A_i\in\calT_2 V$ are alternating tensors whose corresponding skew-symmetric maps $C_i:V\rightarrow V$ satisfy the Clifford commutation relations (\ref{equ1.8}). The Jacobi operators of these tensors (\ref{equ6.2}) are diagonalizable.

In a further hint P. B. Gilkey pointed out that there is a simple possibility to generate examples of Osserman algebraic curvature tensors of the form $\gamma(S)$ or $\alpha(A)$ ($S,A\in\calT_2 V$ symmetric or alternating) if we have a pseudo-Riemannian metric $g$ of indefinite signature $(p,q)$. In particular, we can obtain in this way Osserman algebraic curvature tensors not of the form given in (\ref{osserman1}) or (\ref{equ6.2}) (see Remark \ref{rem6.7}).
\begin{Prop} \label{prop6.5}
Assume that $V$ has a metric $g$ with signature $(p,q)$.
\begin{enumerate}
\item{If $p,q\ge 1$ and $C:V\rightarrow V$ is a symmetric map of $V$ with $C^2 = 0$, then
\begin{eqnarray}
T & := & \gamma(S)\;\;\;,\;\;\;S(x,y)\;:=\;g(Cx,y)\;\;\forall\,x,y\in V \label{equ6.3}
\end{eqnarray}
is an Osserman algebraic curvature tensor.}
\item{If $p,q\ge 2$ and $C:V\rightarrow V$ is a skew-symmetric map of $V$ with $C^2 = 0$, then
\begin{eqnarray}
T & := & \alpha(A)\;\;\;,\;\;\;A(x,y)\;:=\;g(Cx,y)\;\;\forall\,x,y\in V \label{equ6.4}
\end{eqnarray}
is an Osserman algebraic curvature tensor.}
\item{The Jacobi operators $J_T(x)$ of {\rm (\ref{equ6.3})} and {\rm (\ref{equ6.4})} are nilpotent of order {\rm 2}, i.e. $J_T(x)^2 = 0$ for all $x\in V$.}
\end{enumerate}
\end{Prop}
Note that the tensors $S,A\in\calT_2 V$ defined by (\ref{equ6.3}), (\ref{equ6.4}), respectively, are symmetric or skew-symmetric, respectively.

To prove Proposition \ref{prop6.5} we show first the existence of the above nilpotent symmetric or skew-symmetric maps $C$. Furhtermore we show that such maps $C$ do not exist in the positive definite case.
\begin{Lem} \label{lem6.6}
\begin{enumerate}
\item{If $M$ is a symmetric or skew-symmetric $m\times m$-matrix with $M^2 = 0$, then $M = 0$.}
\item{Let $p,q\ge 2$, $m = p+q$, and $F$ be the $m\times m$-diagonal matrix
\begin{eqnarray}
F & := & \mathrm{diag}(\underbrace{1,\ldots1}_{p},\underbrace{-1,\ldots,-1}_{q})
\end{eqnarray}
Then we can find non-vanishing symmetric or skew-symmetric $m\times m$-matrices $S$, $A$, respectively, such that
\begin{eqnarray}
(S\cdot F)^2\;=\;0 &\;\;\;,\;\;\;& (A\cdot F)^2\;=\;0\,.
\end{eqnarray}
The assertion about a symmetric matrix $S$ holds also for $p,q\ge 1$.}
\end{enumerate}
\end{Lem}
\begin{proof}
{\bf Ad (1)}: If a matrix $M = (m_{ij})$ is symmetric or skew-symmetric, we have $M^T = \epsilon\,M$, $\epsilon\in\{1,-1\}$. Then $M^2 = 0$ leads to $\epsilon\,M\cdot M^T = 0$ and
\begin{eqnarray}
0 & = & \mathrm{trace}\,M\cdot M^T \;=\;\sum_{i,j}\,m_{ij}^2\,. \label{equ6.7}
\end{eqnarray}
From (\ref{equ6.7}) we obtain $m_{ij} = 0$ and $M = 0$.

{\bf Ad (2)}: The existence of skew-symmetric matrices $A$ with $(A\cdot F)^2 = 0$ is guaranteed by an example of P. B. Gilkey \cite[p.186]{gilkey5} which was constructed for $p,q\ge 2$.

Let us now consider symmetric matrices under $p,q\ge 1$. If we write $F$ as block matrix
\begin{eqnarray*}
F & = & \left(
\begin{array}{ccc}
\ddots & & \\
 & 
\scriptsize{
\begin{array}{rr}
1 & 0 \\
0 & -1 \\
\end{array}
}
& \\
 & & \ddots\\
\end{array}
\right)
\end{eqnarray*}
and choose a corresponding symmetric block matrix
\begin{eqnarray*}
S & = & \left(
\begin{array}{ccc}
0 & & \\
 & 
\scriptsize{
\begin{array}{cc}
1 & 1 \\
1 & 1 \\
\end{array}
}
& \\
 & & 0\\
\end{array}
\right)\,
\end{eqnarray*}
then we obtain
\begin{eqnarray*}
S\cdot F & = & \left(
\begin{array}{ccc}
0 & & \\
 & 
\scriptsize{
\begin{array}{cc}
1 & -1 \\
1 & -1 \\
\end{array}
}
& \\
 & & 0\\
\end{array}
\right)\,.
\end{eqnarray*}
Obviously, the matrix $S\cdot F$ satisfies $(S\cdot F)^2 = 0$.
\end{proof}
A relation $B(x,y) = g(Cx,y)$ between a linear map $C:V\rightarrow V$ and a tensor $B\in\calT_2 V$ leads to relations
\begin{eqnarray} \label{equ6.8}
B_{ij}\;=\;C_i^k\,g_{kj} & \;\;\;\mathrm{or}\;\;\; & C_i^j\;=\;B_{ik}\,g^{kj}
\end{eqnarray}
for the coordinates of $B$, $C$, and $g$. If we use an orthonormal basis of $V$ then we have $(g^{ij}) = \mathrm{diag}(1,\ldots,1,-1,\ldots,-1)$ and the matrix $(B_{ik}\,g^{kj})$ is a matrix product of the type $B\cdot F$ considered in Lemma \ref{lem6.6}. Thus statement (2) of Lemma \ref{lem6.6} guarantees the existence of symmetric or skew-symmetric maps $C$ with $C^2 = 0$ in the above pseudo-Riemannian settings. Furthermore we see from statement (1) of Lemma \ref{lem6.6} that such maps do not exist if $g$ is positive definite.\\*[0.3cm]
Now we can prove Proposition \ref{prop6.5}.
\begin{proof}
For a symmetric map $C:V\rightarrow V$ the definitions (\ref{jacobi}), (\ref{gammas}) and $S(x,w) = g(Cx,w)$ lead to
\begin{eqnarray*}
g(J_{\gamma(S)}(x)y,w) & = &
{\textstyle \frac{1}{3}}(S(y,w)S(x,x) - S(x,y)S(x,w))
\end{eqnarray*}
and 
\begin{eqnarray}
J_{\gamma(S)}(x)y & = &
{\textstyle \frac{1}{3}}(g(Cx,x)Cy - g(Cy,x)Cx)\,. \label{equ6.9}
\end{eqnarray}
If now $C^2 = 0$, then we obtain
$CJ_{\gamma(S)}(x)y = \frac{1}{3}(g(Cx,x)C^2y - g(Cy,x)C^2x) = 0$ and $J_{\gamma(S)}(x)^2y = \frac{1}{3}(g(Cx,x)CJ_{\gamma(S)}(x)y - g(CJ_{\gamma(S)}(x)y,x)Cx) = 0$, i.e. $J_{\gamma(S)}(x)$ is nilpotent of order 2 for all $x\in V$. But then all eigenvalues of $J_{\gamma(S)}(x)$ are equal to zero, i.e. $\gamma(S)$ is Osserman.

For a skew-symmetric map $C:V\rightarrow V$ we obtain from (\ref{jacobi}), (\ref{alphaa}), $A(x,w) = g(Cx,w)$ and $A(x,x) = 0$
\begin{eqnarray*}
g(J_{\alpha(A)}(x)y,w) & = &
{\textstyle \frac{1}{3}}(2A(y,x)A(x,w) + A(y,x)A(x,w) - A(y,w)A(x,x))\\
 & = & A(y,x)A(x,w)
\end{eqnarray*}
and
\begin{eqnarray*}
J_{\alpha(A)}(x)y & = & g(Cy,x)Cx\,.
\end{eqnarray*}
Again, the condition $C^2 = 0$ leads to $CJ_{\alpha(A)}(x)y = 0$ and
$J_{\alpha(A)}(x)^2y = 0$ for all $x\in V$. Consequently all eigenvalues of $J_{\alpha(A)}(x)$ vanish for all $x\in V$ and $\alpha(A)$ is Osserman.
\end{proof}
\begin{Rem} \label{rem6.7}
If we use a representation of (\ref{osserman1}) or (\ref{equ6.2}) in which we arrange the skew-symmetric maps $C_i$ as in (\ref{equ1.9}) then we will see that the Jacobi operator of (\ref{osserman1}) or (\ref{equ6.2}) has the eigenvalues
\[
0\;\;\;,\;\;\;{\lambda}_0\;\;\;,\;\;\;{\lambda}_0 - 3 {\lambda}_i\,.
\]
on the (pseudo)-unit sphere (see \cite{gilkey3} and \cite[p.191]{gilkey5}). The vanishing of all these eigenvalues leads to $T = 0$ for (\ref{osserman1}) and (\ref{equ6.2}). Proposition \ref{prop6.5}, however, yields examples of non-vanishing Osserman algebraic curvature tensors $\gamma(S)$, $\alpha(A)$ whose Jacobi operators have only the eigenvalues zero. Consequently the examples of Proposition \ref{prop6.5} can not be transformed into a representation (\ref{osserman1}) or (\ref{equ6.2}).

Moreover, the Jacobi operators of the algebraic curvature tensors (\ref{osserman1}) and (\ref{equ6.2}) are diagonalizable (see \cite{gilkey3} and \cite[p.191]{gilkey5}). Since the Jacobi operators $J_{\gamma(S)}$, $J_{\alpha(A)}$ from Proposition \ref{prop6.5} have only the eigenvalue zero, a Jacobi operator $J_{\gamma(S)}$, $J_{\alpha(A)}$ according to Proposition \ref{prop6.5} is equal to zero if it is diagonalizable. Thus every algebraic curvature tensor $\gamma(S)$, $\alpha(A)$ from Proposition \ref{prop6.5}, which has a non-vanishing Jacobi operator, is different from every algebraic curvature tensor with a diagonalizable Jacobi operator.
\end{Rem}
\begin{Rem}
Osserman tensors (\ref{equ6.3}) for metrics $g$ with Lorentzian signature $(1,q)$ are a special case. It is known that an Osserman algebraic curvature tensor has constant sectional curvature if the metric $g$ has Lorentzian signature $(1,q)$, $q\ge 1$ (see \cite{bbg97,gkv97}). The Jacobi operator of an algebraic curvature tensor with constant sectional curvature is diagonalizable. Thus we obtain from Remark \ref{rem6.7} that every algebraic curvature tensor $\gamma(S)$ from Proposition \ref{prop6.5} has a Jacobi operator $J_{\gamma(S)} = 0$ if $p = 1$ or $q = 1$. We verify this fact.
\end{Rem}
\begin{Lem} \label{lem6.9}
Let $F$ be the $m\times m$-diagonal matrix
\begin{eqnarray*}
F & := & \mathrm{diag}(1,-1,\ldots,-1)\,.
\end{eqnarray*}
\begin{enumerate}
\item{If $S\not= 0$ is a symmetric $m\times m$-matrix with $(S\cdot F)^2 = 0$ then there exists an invertible matrix $D$ and a $\lambda\not= 0$ such that
\begin{eqnarray}
D\cdot S\cdot D^{-1} & = & \left(
\begin{array}{cc}
\scriptsize{
\begin{array}{rr}
\lambda & \pm\lambda \\
\pm\lambda & \lambda \\
\end{array}
} & \\
 & 0\\
\end{array}
\right)\;\;\;\;and\;\;\;\;D\cdot F\cdot D^{-1} = F\,. \label{equ6.10a}
\end{eqnarray}
}
\item{A skew-symmetric $m\times m$-matrix $A$ with $(A\cdot F)^2 = 0$ vanishes.}
\end{enumerate}
\end{Lem}
\begin{proof}
{\bf Ad (1):} We write a symmetric $m\times m$-matrix $S$ as block matrix
\begin{eqnarray}
S & = & \left(
\begin{array}{cc}
a & b \\
b^T & \tilde{S}\\
\end{array}
\right)\,, \label{equ6.10}
\end{eqnarray}
where $\tilde{S}$ is a symmetric $(m-1)\times (m-1)$-matrix. Then we can transform $\tilde{S}$ into a diagonal matrix by a conjugation of $S$ by a matrix
\begin{eqnarray*}
D & = & \left(
\begin{array}{cc}
1 & 0\\
0^T & \tilde{D}\\
\end{array}
\right)
\end{eqnarray*}
where $\tilde{D}$ is a suitable orthogonal $(m-1)\times (m-1)$-matrix. Obviously it holds $D\cdot F\cdot D^{-1} = F$ for such a $D$.

From (\ref{equ6.10}) we obtain
\begin{eqnarray*}
S\cdot F\;=\;
\left(
\begin{array}{cc}
a & -b\\
b^T & -\tilde{S}\\
\end{array}
\right)
& \;\;,\;\; &
(S\cdot F)^2 \;=\;
\left(
\begin{array}{cc}
a^2 - b\cdot b^T & -a\,b + b\cdot\tilde{S}\\
a\,b^T - \tilde{S}\cdot b^T & - b^T\cdot b + \tilde{S}^2\\
\end{array}
\right)\,.
\end{eqnarray*}
Assume that $\tilde{S} = 0$. Then the condition $(S\cdot F)^2 = 0$ leads to the relations $a^2 - b\cdot b^T = 0$ and $ab = 0$ from which $S = 0$ follows. However this is a contradiction to $S\not= 0$.

In case of $\tilde{S}\not= 0$ the matrix $\tilde{S}$ has rank 1 because of
$b^T\cdot b - \tilde{S}^2 = 0$. Then $\tilde{S}$ has one and only one non-vanishing eigenvalue $\lambda$. Furthermore we obtain $b\not= 0$ from $b^T\cdot b - \tilde{S}^2 = 0$.

Now the condition $ab^T - \tilde{S}\cdot b^T = 0$ means that $b^T$ is an eigenvector of $\tilde{S}$ which belongs to the eigenvalue $a$. The value $a = 0$ is impossible since then $a^2 - b\cdot b^T = 0$ yields $b = 0$. Thus we obtain $a = \lambda\not= 0$.

If we place $\lambda\not= 0$ in the left upper corner of $\tilde{S}$, then $b$ is an $(m-1)$-tuple $b = (\mu,0,\ldots,0)$ and $a^2 - b\cdot b^T = 0$ yields $\mu = \pm\lambda$.

{\bf Ad (2)}: We write a skew-symmetric $m\times m$-matrix $A$ as block matrix
\begin{eqnarray}
A & = &
\left(
\begin{array}{cc}
0 & b\\
-b^T & \tilde{A}\\
\end{array}
\right)\,, \label{equ6.11}
\end{eqnarray}
where $\tilde{A}$ is a skew-symmetric $(m-1)\times (m-1)$-matrix. From (\ref{equ6.11}) we obtain
\begin{eqnarray*}
A\cdot F\;=\;
\left(
\begin{array}{cc}
0 & -b\\
-b^T & -\tilde{A}\\
\end{array}
\right)
& \;\;,\;\; &
(A\cdot F)^2 \;=\;
\left(
\begin{array}{cc}
b\cdot b^T & b\cdot\tilde{A}\\
\tilde{A}\cdot b^T & b^T\cdot b + \tilde{A}^2\\
\end{array}
\right)\,.
\end{eqnarray*}
Now the condition $(A\cdot F)^2 = 0$ yields $b\cdot b^T = 0$, i.e. $b = 0$. But then $\tilde{A}$ is a skew-symmetric matrix which has to fulfil $\tilde{A}^2 = 0$. From Lemma \ref{lem6.6} we obtain $\tilde{A} = 0$.
\end{proof}
Now let $g$ be a Lorentzian metric with signature $(1,q)$ and $T = \gamma(S)$ be an Osserman algebraic curvature tensor according to Proposition \ref{prop6.5}. Then $S$ can be transformed into a form (\ref{equ6.10a}) and the corresponding symmetric map
\begin{eqnarray*}
C\;:=\;S\cdot F & = & \left(
\begin{array}{cc}
\scriptsize{
\begin{array}{rr}
\lambda & \mp\lambda \\
\pm\lambda & -\lambda \\
\end{array}
} & \\
 & 0\\
\end{array}
\right)
\end{eqnarray*}
has rank 1. Let us consider the Jacobi operator (\ref{equ6.9}). If $Cx = 0$ then $J_{\gamma(S)}(x) = 0$. If $Cx\not= 0$ then the vectors $Cx$ and $Cy$ are proportional since $C$ has a 1-dimensional range. Again we obtain $J_{\gamma(S)}(x) = 0$. Thus $J_{\gamma(S)}(x)\equiv 0$ for every $T = \gamma(S)$ according to Proposition \ref{prop6.5} in the Lorentzian setting. Statement (2) of Lemma \ref{lem6.9} tells us that a non-vanishing Osserman tensor $T = \alpha(A)$ of the type (\ref{equ6.4}) does not exist in the case of a Lorentzian signature.\\*[0.3cm]
\noindent {\bf Acknowledgements.} I would like to thank Prof. P. B. Gilkey for important and helpful discussions and for valuable suggestions for future investigations.
\vspace{0.4cm}
%\bibliography{bflit}
%\bibliographystyle{plain}
%\end{document}

\end{document}